\documentstyle[12pt,a4,pstricks,pst-node,pst-coil]{article}
 \newfont{\Bbb}{msbm10 scaled 1200}
 \newfont{\typewriter}{cmtt10 scaled 1200}
 
 \newcommand{\Z}{\mbox{\Bbb Z}}
 \newcommand{\N}{\mbox{\Bbb N}}
 
 \newcommand{\R}{\mbox{\Bbb R}}
 \newcommand{\C}{\mbox{\Bbb C}}

 \newcommand{\ep}{\vspace{0.5cm}}

\begin{document}

 \title{On Nodal Sets for Dirac and Laplace Operators}
 \author{Christian B\"ar
 \thanks{Partially supported by SFB 256 and by the GADGET program of the EU}}
 \date{March, 1997}
 \maketitle
 
 \begin{abstract}
 \noindent
 We prove that the nodal set (zero set) of a solution of a generalized
 Dirac equation on a Riemannian manifold has codimension 2 at least.
 If the underlying manifold is a surface, then the nodal set is
 discrete.
 We obtain a quick proof of the fact that the nodal set of an
 eigenfunction for the Laplace-Beltrami operator on a Riemannian
 manifold consists of a smooth hypersurface and a singular set of
 lower dimension.
 We also see that the nodal set of a $\Delta$-harmonic differential
 form on a closed manifold has codimension 2 at least; a fact which is not
 true if the manifold is not closed.
 Examples show that all bounds are optimal.
 
  
 \end{abstract}
 
 \section{Introduction}
 The motion of a vibrating membrane $M$ fixed at the boundary is described
 by a function $u: M\times \R \to \R$ satisfying the wave equation
 $$\frac{\partial^2u}{\partial t^2} + \Delta u = 0$$ 
 and Dirichlet boundary conditions $u|_{\partial M} = 0$.
 One can expand $u$ into a series $u(x,t) = \sum_{j=0}^\infty (a_j
 \sin(\sqrt{\lambda_j}t) + b_j\cos(\sqrt{\lambda_j}t))\phi_j(x)$ where 
 $\phi_j$ are the $\Delta$-eigenfunctions on $M$ for the eigenvalue
 $\lambda_j$.
 A ``pure sound'' is given by $u(x,t) =(a_j\sin(\sqrt{\lambda_j}t) + 
 b_j\cos(\sqrt{\lambda_j}t))\phi_j(x)$ for some fixed $j$.
 The zero set of $\phi_j$ describes those points of the membrane which
 do not move during the vibration.
 They can be made visible by putting fine powder on the membrane.
 
 The structure of the zero set of eigenfunctions of the
 Laplace-Beltrami operator on surfaces is well understood.
 They consist of smooth arcs, called nodal lines, and isolated
 singular points where these arcs meet.
 One knows that the arcs meeting at a singular point form an
 equiangular configuration, see \cite{Al} or \cite[Satz 1]{Br} for a 
 proof.
 One also has lower and upper bounds for the length of the nodal lines
 \cite[Satz 2]{Br}, \cite[Theorem 1.2, Corollary 1.3]{DF3},
 \cite[Theorem 4.2]{Dg}.
 
 Such a precise understanding of nodal sets seems difficult in higher
 dimensions.
 But one has the following regularity result.
 
 \ep
 
 {\em If $\phi$ is an eigenfunction of the Laplace-Beltrami operator on
 an $n$-dimensional Riemannian manifold, then the nodal set of $\phi$
 consists of a smooth hypersurface and a singular part of dimension
 $\le n-2$.}
 
 \ep
 
 Essentially this fact has been stated as Theorem 2.2 in \cite{Ch}.
 Once it is established the standard proof of Courant's nodal domain
 theorem in dimension~2 carries over to higher dimensions.
 A nodal domain is a connected component of the complement of
 the nodal set.
 Courant's nodal domain theorem states that the number of nodal
 domains of the $i^{th}$ eigenfunction of $\Delta$ is less
 than or equal to $i$.
 The point is that one has to use a Green formula over a nodal
 domain and this requires some regularity of its boundary. 
 
 It was pointed out by Y.\ Colin de Verdi\`ere that the proof of
 Theorem 2.2 in \cite{Ch} has a serious gap, see \cite[App.\ E]{BM}.
 Therefore B\'{e}rard and Meyer modified the proof of Courant's nodal
 domain theorem.
 They approximate the nodal domains by regular domains \cite[App. D]{BM}.
 
 The above regularity statement on Laplace-eigenfunctions has (to our
 knowledge) first been proved by R.\ Hardt and L.\ Simon in
 \cite[Theorem 1.10]{HS}, see also \cite{CF} for a special case.
 We will obtain a quick proof of this fact once our theorem on the
 nodal set of solutions of Dirac equations is established (Corollary 2).
 
 One has estimates for the $(n-1)$-dimensional Hausdorff measure
 of the nodal set, see \cite[Theorem 4.2]{Do}, \cite[Theorem 1.2]{DF1},
 \cite{DF2}, and \cite[Theorem 5.3]{HS}, some of which only work if
 the Riemannian metric is assumed to be real analytic.
 There are also estimates for the volume of the nodal domains, see
 \cite[Theorem B]{Lu}, \cite[Theorem 2]{CM}.
 
 Not much is known about the topology of the singular set.
 The problem is that in dimension $n \ge 3$ the nodal set need not
 be locally homeomorphic to its tangent cone.
 This precisely was the problem in \cite{Ch}.
 For a structural result in dimension $n=3$ see \cite[Theorem 1.2]{Cn}.
 
 The purpose of this paper is to study the nodal set of solutions of
 certain systems of elliptic linear partial differential equations of first
 order, {\em generalized Dirac equations}.
 Precise definitions will be given in the next section.
 
 To get some feeling of what to expect let us first look at the
 trivial one-dimensional case.
 A Laplace equation is then nothing but a second order linear
 ordinary differential equation.
 The standard theory tells us that zeros of solutions must be isolated.
 A Dirac equation becomes a first order linear ordinary differential 
 equation and we see that nontrivial solutions have no zeros at all.
 
 Another interesting test case is provided by holomorphic functions on
 Riemann surfaces.
 The Cauchy-Riemann equations are special generalized Dirac
 equations.
 As is well-known, zeros of solutions must form a discrete set.
 In higher dimensions, a holomorphic function defines a complex
 subvariety having real codimension 2.
 
 This, as well as other special cases, leads us towards the conjecture
 that the nodal set of a solution of a generalized Dirac equation
 on an $n$-dimensional manifold has dimension $\le n-2$.
 
 The main result of this paper says that this conjecture is indeed
 true.
 Examples show that this bound is optimal.
 In the two-dimensional case we have a similar unique-continuation
 theorem as we have for holomorphic functions (Corollary~3).
 If the nodal set of a solution of a generalized Dirac equation on
 a connected surface
 has a cumulation point, then this solution must vanish identically.
 
 There is a physical interpretation in quantum mechanics similar to 
 the vibrating membrane for the Laplace equation.
 The ``wave function'' of a fermion (e.g.\ an electron) is given
 by a spinor field, in the simplest case by a function $\psi : \R^3 
 \times \R \to \C^2$, satisfying the equation 
 $$i\frac{\partial}{\partial t}\psi + (D+h)\psi = 0$$
 where $D$ is the Dirac operator on $\R^3$ and $h$ is a potential.
 The particle is in ``pure state'' if $\psi$ has the form
 $\psi(x,t) = e^{i\lambda t}\Psi(x)$ where $\Psi$ is an eigenfunction
 of $D+h$ for the eigenvalue $\lambda$.
 The scalar function $|\psi(x,t)|^2 = |\Psi(x)|^2$ can be interpreted
 as the probability measure for the particle to be found at the point $x$.
 Hence the nodal set of $\Psi$ is the set of points where this 
 probability measure is zero.
 Our main result then says that this ``exclusive set'' is not very big, it 
 is at most one-dimensional.
 
 We believe that the nodal set of solutions of Dirac equations
 very often carries
 important information about the underlying manifold.
 In the special case of holomorphic functions on complex manifolds
 this is classical.
 Recently, Taubes has given an impressive example in the theory
 of 4-dimensional symplectic manifolds.
 Starting from solutions $\psi_r$ of a one-parameter family
 of Dirac equations (deformed first Seiberg-Witten equation),
 $r \in \R$, he constructs pseudoholomorphic curves.
 Philosophically, these curves are given by the zero locus of 
 $\psi_\infty$ (which is 2-dimensional!).
 See \cite{T1}, \cite{T2}, and \cite{Ko} for a very readable
 survey. 
 
 It is nice that information about solutions of Dirac equations 
 yields also information about solutions of Laplace equations.
 It was mentioned earlier that we will obtain a quick proof of the fact
 that the nodal set of a $\Delta$-eigenfunction is the union of a 
 smooth hypersurface and a singular part of dimension $\le n-2$.
 
 We will also see that the nodal set of a $\Delta$-harmonic
 differential form on a closed manifold has codimension 2 at least
 (Corollary 1).
 This is surprising since it is not true if the underlying manifold
 is not closed nor does it hold for other $\Delta$-eigenforms even if 
 the manifold is closed.
 
 The paper is organized as follows.
 In the next section we give precise definitions and we collect a
 few well-known facts about Dirac operators for later use.
 We then formulate the main result.
 
 In the third section we give the most important examples for 
 generalized Dirac operators and for some of these we draw
 conclusions from our main theorem.
 
 In section 4 we give the proof.
 We employ tools similar to those used in complex algebraic geometry
 when one studies the topological structure of complex algebraic
 varieties.
 We use an analog of Weierstrass' preparation theorem for 
 differentiable functions due to Malgrange to write a solution of 
 a Dirac equation in a certain normal form.
 This theorem is important e.g.\ in catastrophe theory.
 To insure that we can apply this theorem we have to use Aronszajn's
 unique continuation theorem.
 We easily conclude that the nodal set has dimension not bigger than
 $n-1$.
 
 The difficult part is to show that simultaneous vanishing of
 several components of our section must reduce the dimension of its
 zero set once more.
 This requires a careful investigation of certain resultants
 and this is where we really use the Dirac equation.

 \section{Statement of Result}
 Let $M$ be an $n$-dimensional Riemannian manifold,
 let $S$ be a Riemannian or Hermitian vector bundle over $M$
 on which the Clifford bundle $Cl(TM)$ acts from the left.
 This means that at every point $p \in M$ there is a linear
 map $T_pM \otimes S_p \to S_p, v \otimes s \to v \cdot s$, 
 satisfying the relations 
 $$
 v \cdot w \cdot s + w \cdot v \cdot s = -2 \langle v,w\rangle s
 $$
 and 
 $$\langle v \cdot s_1 , s_2 \rangle = - \langle s_1 , 
 v \cdot s_2 \rangle .
 $$
 Moreover, let $\nabla$ be a metric connection on $S$ satisfying 
 the Leibniz rule for Clifford multiplication
 $$
 \nabla (v \cdot s) = (\nabla v)\cdot s + v \cdot \nabla s.
 $$
 Such an $S$ is called a {\em Dirac bundle}.
 All geometric data such as metrics and Clifford multiplication are
 assumed to be $C^\infty$-smooth.
 
 The {\em (generalized) Dirac operator} acts on sections of $S$ and 
 is defined by 
 $$
 Ds = \sum_{i=1}^n e_i \cdot \nabla_{e_i} s
 $$ 
 where $e_1 , \ldots , e_n$ denote a local orthonormal frame of $TM$.
 The Dirac operator is easily seen to be independent of the choice
 of orthonormal frame.
 It is a first order formally self-adjoint elliptic differential
 operator, see \cite{BGV} or \cite{LM} for details.
 
 Three general facts about Dirac operators will be used
 later on and should be stated here for completeness.
 The proofs are simple computations.
 
 {\bf 1.} If $s$ is a differentiable section of $S$ and $f$ a differentiable
 function on $M$, then the following {\em Leibniz rule} holds 
 \cite[p.\ 116, Lemma 5.5]{LM}:
 \begin{equation}
 D(f\cdot s) = f\cdot Ds + \nabla f\cdot s
 \label{DPR} 
 \end{equation}
 where $\nabla f$ is the gradient of $f$.
 
 {\bf 2.} The connection $\nabla$ maps sections of $S$ into those of 
 $T^\ast M\otimes S$.
 Denote the formal adjoint of this operator by $\nabla^\ast$.
 Then we have the following {\em Weitzenb\"ock formula} 
 \cite[p.\ 155, Theorem 8.2]{LM}:
 \begin{equation}
 D^2 = \nabla^\ast\nabla + \Re
 \label{WBF} 
 \end{equation}
 where $\Re$ is an endomorphism field given by curvature.
 Operators of the form $\nabla^\ast\nabla + 
 \mbox{endomorphism field}$ are called {\em generalized Laplacians}.
 
 {\bf 3.} If $M$ has smooth boundary $\partial M$ with outer unit
 normal field $\nu$, if $s_1$ and $s_2$ are compactly
 supported $C^1$-sections of $S$, then the following 
 {\em Green formula} holds \cite[p.\ 115, eq.\ (5.7)]{LM}:
 \begin{equation}
 (Ds_1,s_2)_{L^2(M)} - (s_1,Ds_2)_{L^2(M)} =
 \int_{\partial M} \langle \nu\cdot s_1,s_2 \rangle .
 \label{DGF} 
 \end{equation}
 
 \ep
 
 If $s$ is an eigensection of $D$, 
 i.e.\ $Ds = \lambda s, \lambda \in \R$, or, more generally, 
 $s$ satisfies $(D + h)s = 0$ for
 some endomorphism field $h$ of $S$, then the zero locus of $s$, 
 $\{ x\in M\ |\ s(x) = 0 \}$,
 is called the {\em nodal set} of $s$.
 The main purpose of this paper is to study the structure of such
 nodal sets.
 \ep
 
 {\bf Main Theorem.}
 {\em Let $M$ be a connected $n$-dimensional Riemannian manifold 
 with Dirac bundle $S$ and generalized Dirac operator $D$.
 Let $h$ be a smooth endomorphism field for $S$ and
 let $s \not\equiv 0$ be a solution of 
 $$
 (D + h)s = 0.
 $$
 Then the nodal set of $s$ is a countably $(n-2)$-rectifiable set
 and thus has Hausdorff dimension $n-2$ at most.
 
 If $n=2$, then the nodal set of $s$ is a discrete subset
 of $M$.
 }
 \ep
 
 Recall that a subset of an $n$-dimensional Riemannian manifold $M$
 is called {\em countably k-rectifiable} if it can be written as
 a countable union of sets of the form $\Phi(X)$ where $X \subset
 \R^k$ is bounded and $\Phi : X \to M$ is a Lipschitz map.
 The proof will be given in the fourth section.

 \section{Examples and Consequences}
 Let us look at the most important examples and draw
 some conclusions.
 
 \ep
 
 {\bf Example 1.}
 The Clifford algebra bundle $Cl(TM)$ acts on itself by 
 Clifford multiplication.
 As a vector bundle $Cl(TM)$ can be canonically identified
 with the exterior form bundle $\Lambda^\ast (TM)$.
 We thus obtain a real Dirac bundle $S = \Lambda^\ast (TM)$
 with Levi-Civita connection $\nabla$.
 The Dirac operator is $D = d + \delta$ where $d$ denotes
 exterior differentiation and $\delta$ codifferentiation.
 
 \ep
 
 {\bf Example 2.}
 Let $M$ be a spin manifold.
 Then there exists the spinor bundle $S$, a complex Dirac
 bundle of rank $2^{[n/2]}$.
 The connection $\nabla$ is induced by the Levi-Civita connection.
 The corresponding Dirac operator $D$ is the classical Dirac
 operator.
 
 \ep
 
 {\bf Example 3.}
 Let $M$ be a spin$^c$ manifold.
 Again, there exists the complex spinor bundle $S$ of rank $2^{[n/2]}$.
 The connection $\nabla$ depends on the Levi-Civita connection
 and the choice of a connection on a certain U(1)-bundle,
 the {\em determinant bundle}.
 
 \ep
 
 {\bf Example 4.}
 Let $M$ be an almost complex manifold.
 Then $M$ is canonically spin$^c$ and the spinor bundle of example 3 
 can be identified with the bundle of mixed $(0,p)$-forms, $S =
 \Lambda^{0,\ast}(TM\otimes \C)$.
 The Dirac operator is given by $D = \sqrt{2} \cdot (\bar{\partial} +
 \bar{\partial}^\ast)+h$ where $h$ is a zero order term
 which vanishes if $M$ is K\"ahler.
 
 \ep
 
 {\bf Example 5.}
 Let $S$ be a Dirac bundle over $M$, and let $E$ be a Riemannian
 or Hermitian bundle over $M$ with metric connection.
 Then $S \otimes E$ canonically becomes a Dirac bundle and the
 corresponding Dirac operator is called {\em twisted Dirac operator
 with coefficients in E}.
 
 \ep
 
 Let us see what the theorem tells us when applied to the example
 of differential forms.
 First of all, we see that the bound in the theorem is optimal.
 Namely, let $F$ be a surface of higher genus, let $\omega_1$ be a 
 nontrivial closed and coclosed 1-form on $F$.
 In other words, $\omega_1$ is a solution of $(d + \delta)\omega_1 = 0$.
 The nodal set of $\omega_1$ consists of isolated points and it is
 nonempty since the Euler number of $F$ is nonzero.
 Let $T^{n-2}$ be a flat $(n-2)$-torus and let $\omega_2$ be a parallel
 $k$-form on $T^{n-2}$.
 Put $M = F \times T^{n-2}$ and denote the projections onto the
 factors by $\pi_1 : M \to F$ and $\pi_2 : M \to T^{n-2}$.
 Then $\omega = \pi_1^\ast\omega_1 \wedge \pi_2^\ast\omega_2$ is a
 closed 
 and coclosed $(k+1)$-form on $M$.
 Its nodal set is a disjoint union of copies of $T^{n-2}$ and therefore
 has codimension 2 in $M$.
 
 \ep
 
 {\bf Corollary 1.}
 {\em Let $M$ be a compact connected Riemannian manifold
 without boundary.
 Let $\Delta = d\delta + \delta d$ be the Laplace-Beltrami operator
 acting on $k$-forms.
 Let $\omega \not\equiv 0$ be a harmonic $k$-form, i.e.\ $\Delta\omega
 = 0$.
 
 Then the nodal set of $\omega$ is a countably $(n-2)$-rectifiable set
 and thus has Hausdorff dimension $n-2$ at most.
 
 If $n=2$, then the nodal set of $\omega$ is a discrete
 subset of $M$.
 }
 
 \ep
 
 P{\footnotesize ROOF}.
 Taking $L^2$-products and using the Green formula (\ref{DGF}) yields
 \begin{eqnarray*}
 0 & = & (\Delta \omega ,\omega )_{L^2(M)} \\
 & = & ((d + \delta)^2\omega ,\omega )_{L^2(M)} \\
 & = & ((d + \delta)\omega ,(d + \delta)\omega )_{L^2(M)}
 .
 \end{eqnarray*}
 We conclude
 $$
 (d + \delta)\omega \equiv 0.
 $$
 Now the theorem tells us that the nodal set of $\omega$ is a countably
 $(n-2)$-rectifiable set.
 $\Box$
 
 \ep
 
 It is remarkable that Corollary 1 fails if we drop the assumption that
 $M$ be compact.
 For example, $\omega = x_1 dx_1 \wedge \ldots \wedge dx_k$ is a
 harmonic $k$-form on $\R^n$ whose nodal set has codimension 1.
 This means that despite the local nature of the statement of Corollary
 1 it can not be proved by purely local methods.
 
 One can extend Corollary 1 to complete noncompact manifolds
 by imposing suitable decay conditions on the form $\omega$.
 This will make the additional boundary term 
 $\int_{B(R)} \langle \nu\cdot(d + \delta)\omega ,\omega\rangle$
 tend to zero as $R\to\infty$ where $B(R)$ denotes the distance
 ball of radius $R$ around some fixed point.
 Demanding $\omega$ to be in the Sobolev space 
 $H^{\frac{1}{2},2}(M)$ should be sufficient.
 
 Corollary 1 also fails if we replace harmonic forms by other
 eigenforms of the Laplace-Beltrami operator even if the manifold
 is compact.
 For example, $\omega = \sin(2\pi m x_1) dx_1 \wedge \ldots \wedge dx_k$,
 $m \in \Z$, is a $\Delta$-eigenform on the torus $T^n = \R^n/\Z^n$
 whose nodal set has codimension 1.
 
 \ep
 
 As another application of the main theorem we obtain a quick proof
 of the following theorem first proven 1989 by R.\ Hardt and L.\ 
 Simon \cite[Theorem 1.10]{HS}.
 
 \ep
 
 {\bf Corollary 2.}
 {\em Let $M$ be an $n$-dimensional connected Riemannian manifold.
 Let $f$ be a nontrivial $\Delta$-eigenfunction on $M$, i.e.
 $$
 \Delta f = \lambda f
 $$
 for some $\lambda > 0$.
 
 Then the nodal set of $f$ is a disjoint union
 $N_{reg} \cup N_{sing}$ where $N_{reg}$ is a smooth 
 hypersurface of $M$ and $N_{sing}$ is a countably $(n-2)$-rectifiable set
 and thus has Hausdorff dimension $n-2$ at most.}
 
 \ep
 
 P{\footnotesize ROOF}.
 Put $N_{reg} = \{ x \in M\ |\ f(x) = 0, df(x) \not= 0 \}$ and 
 $N_{sing} = \{ x \in M\ |\ f(x) = 0, df(x) = 0 \}$.
 Then the nodal set of $f$ is given by $N_{reg} \cup N_{sing}$. 
 By the implicit function theorem $N_{reg}$ is a smooth hypersurface.
 
 The other set $N_{sing}$ is the zero locus of the mixed differential
 form $\omega = \sqrt{\lambda}f + df$.
 Now $\omega$ is an eigenform for the generalized Dirac operator 
 $D = d + \delta$,
 $$
 (d + \delta)\omega = \sqrt{\lambda}\omega .
 $$
 The main theorem says $N_{sing}$ is a countably $(n-2)$-rectifiable set. 
 $\Box$
 
 \ep
 
 {\bf Remark.}
 If we knew that Corollary 2 is true for generalized
 Laplacians, then the main theorem could be derived from it, 
 at least for eigensections $s$ of a Dirac operator, as follows.
 
 Let $Ds = \lambda s, \lambda \in \R$.
 The Weitzenb\"ock formula (\ref{WBF}) yields
 $$
 \nabla^\ast\nabla s + (\Re - \lambda^2)s = 
 (D^2 - \lambda^2)s = 0.
 $$
 Hence we could conclude that the nodal set of $s$ is of the
 form $N_{reg} \cup N_{sing}$ where $N_{reg}$ is a smooth
 hypersurface and $N_{sing}$ has at least codimension 2.
 It would remain to show $N_{reg} = \emptyset$.
 
 Assume there is $x_0 \in N_{reg}$.
 Choose a small ball $B$ around $x_0$ which is cut by $N_{reg}$
 into two pieces $B_1$ and $B_2$.
 Define 
 $$
 \tilde{s}(x) = \left\{
 \begin{array}{l}
 s(x), \mbox{ if }x \in B_1,\\
 0, \mbox{ if }x \in B_2 .
 \end{array}
 \right.
 $$
 Then $\tilde{s}$ is a continuous section of $S$ over $B$.
 Let $\phi$ be a smooth test section of $S$ with compact
 support contained in $B$.
 Let $\nu$ be the unit normal field of $N_{reg}\cap B$
 pointing into $B_2$.
 
 \begin{center}
 $$
 \pspicture(0,2)(14,8)

 \pscustom[linewidth=2pt]{
 \pscircle(7,5){2}
 \fill[fillstyle=solid,fillcolor=lightgray]
 }
 
 \pscustom[linewidth=0pt]{
 \gsave
 \pscurve(10,4)(10,7)(8,8)(7,7)(7,5)(7,3)(6,2)(6,1)(10,2)
 \fill[fillstyle=solid,fillcolor=white]
 \grestore
 }
 \pscurve[linewidth=2pt](8,8)(7,7)(7,5)(7,3)(6,2)
 \pscircle[linewidth=2pt](7,5){2}
 
 \qdisk(7,5){3pt}
 \uput[180](7,5){$x_0$}
 
 \put(6.1,5.9){$B_1$}
 \put(5,3.2){$B$}
 \put(8.1,4.9){$B_2$}
 \put(7,8){$N_{reg}$}
 
 \psline[linewidth=1pt]{->}(7.1,3.5)(7.6,3.5)
 \psline[linewidth=1pt]{->}(7.05,4.25)(7.55,4.3)
 \psline[linewidth=1pt]{->}(7,5)(7.5,5.05)
 \psline[linewidth=1pt]{->}(6.95,5.75)(7.45,5.8)
 \psline[linewidth=1pt]{->}(6.9,6.6)(7.4,6.55)
 \put(7.6,5.8){$\nu$}
 
 \endpspicture
 $$
 {\bf Fig.\ 1}
 \end{center}
 
 Then by (\ref{DGF})
 \begin{eqnarray*}
 (\tilde{s}, D\phi)_{L^2(B)} & = &
 (s,D\phi)_{L^2(B_1)} \\
 &=& (Ds, \phi)_{L^2(B_1)} + \int_{B\cap N_{reg}}
 \langle s, \nu\cdot\phi \rangle \\
 &=& (\lambda s, \phi)_{L^2(B_1)} + 0\\
 &=& (\lambda \tilde{s}, \phi)_{L^2(B)}.
 \end{eqnarray*}
 Hence $D\tilde{s} = \lambda \tilde{s}$ holds in $B$ in the sense of
 distributions.
 By elliptic regularity theory $\tilde{s}$ is smooth.
 Since $\tilde{s}$ vanishes identically on $B_2$ we know
 from Aronszajn's unique continuation theorem (see next
 section) that $\tilde{s} \equiv 0$ on $B$.
 Hence $s \equiv 0$ on $B_1$.
 Applying Aronszajn's theorem once more we conclude
 $s \equiv 0$ on $M$. 
 $\Box$
 
 \ep
 
 Unfortunately, to our knowledge Corollary 2 is not established 
 for generalized Laplacians. 
 But, using the methods of the next section, it is easy to 
 see that nodal sets for generalized Laplacians
 have Hausdorff dimension $n-1$ at most.
 
 Note that the nodal set of a solution of a general linear elliptic
 system of second order can be very irregular.
 For example, it is not hard to show the following.
 
 {\em Let $A \subset \R^{n-1}$ be any closed subset.
 Then there is a linear elliptic differential operator of second
 order, $P$, acting on functions $u : \R^n \to \R^2$ and there is a solution
 $u$ of $Pu = 0$ such that $u^{-1}(0) = A \times \{0\} \subset \R^n$.}
 
 \ep
 
 To conclude this section let us emphasize once more the 
 two-dimensional case.
 We have the following generalization of the well-known uniqueness
 theorem for holomorphic functions.
 
 \ep
 
 {\bf Corollary 3.}
 {\em Let $M$ be a two-dimensional connected Riemannian manifold.
 If the zero set of a solution of a generalized Dirac equation on $M$
 has a cumulation point, then this solution must vanish identically.}
 $\Box$
 
 \ep
 
 This corollary has been proven for many special cases.
 The Dirac equation on a surface can be written as a generalized
 Cauchy-Riemann equation.
 If for example the real dimension of the Dirac bundle is 2, then
 the theory of ``generalized analytic functions'' applies, see
 \cite{Ve}, also compare \cite{Ca}.

 \section{Proof of the Main Theorem}
 
 This section is devoted to the proof of the main theorem.
 There are two important ingredients to the proof which we
 state first.
 We have
 
 \ep
 
 {\bf Aronszajn's Unique Continuation Theorem.}
 {\em Let $M$ be a connected Riemannian manifold.
 Let $L$ be an operator of the form $L = \nabla^\ast\nabla + L_1 + L_0$
 acting on sections of a vector bundle $S$ over $M$ where $L_1$ and
 $L_0$ are differential operators of first and zero$^{th}$ order respectively.
 Let $s$ be a solution of}
 $$
 Ls = 0.
 $$
 {\em If $s$ vanishes at some point of infinite order, i.e.\ if all
 derivatives vanish at that point, then $s \equiv 0$.} 
 
 \ep
 
 For a proof see \cite[Theorem on p.\ 235 and Remark 3 on p.\ 248]{Ar}.
 The other essential ingredient is a version of Weierstrass' 
 preparation theorem for differentiable functions 
 \cite[Chapter V]{Ma}, \cite[6.3]{BL}.
 
 \ep

 {\bf Malgrange's Preparation Theorem} (Special case).
 {\em Let $U \subset \R^n$ be an open neighborhood of $0$, let $f : U
 \to \R$ be a $C^\infty$-function vanishing of $k^{th}$ order at $0$
 but not of $(k+1)^{st}$ order, $k \in \N$.
 Then, after possibly shrinking $U$ to a smaller neighborhood and
 applying a linear coordinate transformation to $\R^n$, there exist
 $C^\infty$-functions $v : U \to \R$ and $u_j : U \cap (\{0\}\times
 \R^{n-1}) \to \R$ such that
 $$
 f(x) = v(x) \cdot \left( x_1^k +
   \sum_{j=0}^{k-1}u_j(x')x_1^j\right)
 $$ 
 for all $x=(x_1,x') \in U$ where $v(x) \not= 0$ for all $x \in U$
 and $u_j$ vanishes of order $k-j$ at $0 \in \R^{n-1}$.}
 
 \ep
 
 {\bf Remark.}
 By Taylor's theorem we can write $f = \hat{f} + \psi$ where $\hat{f}$
 is a homogeneous polynomial of degree $k$ and $\psi$ vanishes of
 order $k+1$ at 0.
 The linear coordinate transformation in the Preparation Theorem
 must be such that a vector $w\in\R^n$ for which $\hat{f}(w)\not=0$
 is transformed into $(1,0,\ldots,0)$.
 
 \ep
 
 P{\footnotesize ROOF OF} M{\footnotesize AIN} T{\footnotesize HEOREM}.
 To prove the main theorem let $s$ be a section of the Dirac bundle
 satisfying
 \begin{equation}
 (D+h)s=0.
 \label{de}  
 \end{equation}
 We assume that the Dirac bundle is real, in the complex case we
 simply forget the complex structure.
 
 Let $p\in M$ be a point of its nodal set, i.e.\ $s(p)=0$.
 Applying $D$ to (\ref{de}) and using the Weitzenb\"ock formula we get
 (\ref{WBF})
 \begin{equation}
 (\nabla^\ast\nabla + D\circ h + \Re)s = 0.
 \end{equation}
 By Aronszajn's unique continuation theorem $s$ cannot vanish of 
 infinite order at $p$.
 Say $s$ vanishes at $p$ of order $k$ but not of order $k+1$.
 
 We choose normal coordinates around $p$ and trivialize the Dirac
 bundle.
 Then $s$ corresponds to a vector valued function $(s_1,\ldots,s_r)$
 defined in a neighborhood of $0$.
 Here $r$ is the (real) rank of the Dirac bundle $S$.
 All component functions vanish of order $k$ at $p$ and at least
 one of them does not vanish of order $k+1$.
 Other components could a-priori vanish of higher order but by
 choosing the trivialization appropriately we can assume that this is
 not the case.
 To see this, note that trivializing the bundle amounts to exhibiting
 linearly independent linear functionals on the bundle.
 Pick one linear functional $l_1$ such that $s_1 = l_1\circ s$ does not
 vanish of order $k+1$ at $0$.
 Choose the other $r-1$ functionals linearly independent but close
 to $l_1$.
 
 Moreover, if the Dirac bundle is trivialized in this way, then 
 there is a direction in which the $k^{th}$ derivative of all
 components $s_1 , \ldots ,s_r$ does not vanish.
 Hence we can use Malgrange's preparation theorem with the same
 linear coordinate transformation (the same $x_1$-direction) 
 for all the components.
 Therefore we can write
 $$
 s_m(x) = v_m(x) \cdot \left( x_1^k +
   \sum_{j=0}^{k-1}u_{m,j}(x')x_1^j\right)
 $$ 
 where $v_m$ are nonvanishing and $u_{m,j}$ vanish of order $k-j$ at
 $0$, $m = 1,\ldots,r$.
 
 We see already that the nodal set of $s$
 near $p$ is countably $(n-1)$-rectifiable.
 Namely, for any $m$ the nodal set is contained in the set
 $$
 N_m = \left\{x=(x_1,x')\ \left|\ x_1^k +
   \sum_{j=0}^{k-1}u_{m,j}(x')x_1^j = 0 \right.\right\}.
 $$
 That this set $N_m$ is countably $(n-1)$-rectifiable follows easily
 from the following
 
 \ep
 
 {\bf Fact.}
 $\R^k$ can be written as a disjoint union of countably many bounded
 subsets $A_\nu$ such that the number of pairwise distinct real roots
 of the polynomial $P_u(t) = t^k + \sum_{j=0}^{k-1}u_{j}t^j$ is constant
 for $u=(u_0,\ldots,u_{k-1}) \in A_\nu$ and the real roots (ordered my 
 magnitude) are Lipschitz-functions on $A_\nu$.
 
 
 \begin{center}
 $$
 \pspicture(0,1)(14,7)
 
 \pscurve[linewidth=2pt](5.5,2)(6.1,3)(7,4)(8,4.8)(10,5.6)
 \pscurve[linewidth=2pt](7,4)(8,4.4)(9,5.5)(9.5,6.4)
 \pscurve[linewidth=2pt](7,4)(8,3.6)(9,2.5)(9.5,1.6)
 
 \put(10.2,3.4){$x'\in \R^{n-1}$}
 \put(5.5,6.5){$x_1\in \R$}
 \put(8.8,4.6){$N_m$}
 
 \psline[linewidth=1pt]{->}(2,4)(12,4)
 \psline[linewidth=1pt]{->}(7,1)(7,7)
 
 \endpspicture
 $$
 {\bf Fig.\ 2}
 \end{center}
 
 \ep
 
 Recall that two polynomials $F = \sum_{j=0}^k a_jt^j$ and $G = \sum_{j=0}^k b_jt^j$
 have a common root if and only if the {\em resultant} $R_{F,G}$
 vanishes.
 The resultant is a weighted homogeneous polynomial of degree $k^2$ in 
 the coefficients $a_j$ and $b_j$ where $a_j$ and $b_j$ have weight
 $k-j$, see \cite[Section 4]{BK}.
 
 More generally, $r$ polynomials $P_m = \sum_{j=0}^{k}u_{m,j}t^j,
 m=1,\ldots,r$ have a common root if and only if any two linear
 combinations $F = \sum_{m=1}^r \alpha_m P_m$ and $G = \sum_{m=1}^r 
 \beta_m P_m$ have vanishing resultant $R_{F,G}$.
 
 We will show that there are two linear combinations $F = \sum_{m=1}^r
 \alpha_m v_m(0) P_m$ and $G = \sum_{m=1}^r \beta_m v_m(0) P_m$ of our polynomials
 $P_m = t^k + \sum_{j=0}^{k-1}u_{m,j}(x')t^j$ such that the resultant
 $R_{F,G}$ regarded as a function in $x'$ does not vanish of infinite
 order at $x'=0$.
 Then, by applying once more Malgrange's preparation theorem and the
 fact on the real roots of polynomials mentioned above, we see
 that the zero locus of $R_{F,G}$ is a countably $(n-2)$-rectifiable
 subset of $\R^{n-1}$ and that $N = \cap_{m=1}^r N_m$ is contained in 
 a countably $(n-2)$-rectifiable subset of $\R^n$.
 Moreover, if $n=2$, we see that $x=0$ is an isolated zero of
 $s$ and the theorem is proved.
 
 To find two such linear combinations $F$ and $G$ we look at the
 Taylor expansion $v_m(x') = v_m(0) + \mbox{first order terms}$,
 $u_{m,j}(x') = \hat{u}_{m,j}(x') + \mbox{higher order terms}$.
 Here $\hat{u}_{m,j}(x')$ is a homogeneous polynomial of degree $k-j$.
 The Taylor expansion of $s_m$ is then given by
 $$
 s_m(x) = v_m(0) \cdot \left( x_1^k +
   \sum_{j=0}^{k-1}\hat{u}_{m,j}(x')x_1^j\right)
   + \mbox{ higher order terms}.
 $$
 Looking at the lowest order term of the left hand side of equation
 (\ref{de}) we get
 \begin{equation}
 \hat{D}w = 0
 \label{hde} 
 \end{equation}
 where $w(x) = (w_1(x),\ldots,w_r(x)), w_m(x) = v_m(0) \cdot \left( x_1^k +
 \sum_{j=0}^{k-1}\hat{u}_{m,j}(x')x_1^j\right)$ and $\hat{D} =
 \sum_{i=1}^n \gamma_i \partial/\partial x_i$ is the Dirac
 operator on $T_pM = \R^n$.
 Here $\gamma_i$ are generalized Pauli matrices satisfying the
 relations $\gamma_i\gamma_j + \gamma_j\gamma_i = -2\delta_{ij}$.
 
 Putting $y_j(x') = 
 (v_1(0)\hat{u}_{1,j}(x'),\ldots,v_r(0)\hat{u}_{r,j}(x'))$ for 
 $j=0,\ldots,k-1$ and $y_k(x') = (v_1(0),\ldots,v_r(0))$ we have 
 \begin{equation}
 w(x) = \sum_{j=0}^{k}y_j(x')x_1^j.
 \label{ffw} 
 \end{equation}
 
 Write $D_1 = \sum_{j=2}^n\gamma_1\gamma_j\partial/\partial x_j$.
 Then $D_1$ is a generalized Dirac operator on $\R^{n-1}$.
 Plugging (\ref{ffw}) into (\ref{hde}) and comparing
 coefficients of $x_1^j$ yields
 \begin{eqnarray*}
 D_1y_{j} & = & (j+1)\cdot y_{j+1}, 
 \hspace{1cm} j = 0,\ldots,k-1 \\
 D_1y_k & = & 0.
 \end{eqnarray*}
 Hence
 \begin{equation}
 y_j = \frac{1}{j!}D_1^jy_0,
 \hspace{1cm} j = 0,\ldots,k.
 \label{rf} 
 \end{equation}
 
 Given two linear combinations $\hat{F} = \sum_{m=1}^r\alpha_m w_m$
 and $\hat{G} = \sum_{m=1}^r\beta_m w_m$ the resultant
 $R_{\hat{F},\hat{G}}$ is precisely the lowest order term of the 
 resultant of the corresponding linear combinations of the $P_m$'s, $F
 = \sum_{m=1}^r\alpha_mv_m(0)P_m$, $G = \sum_{m=1}^r\beta_mv_m(0)P_m$,
 \begin{equation}
 R_{F,G}(x') = R_{\hat{F},\hat{G}}(x') + \mbox{ higher order terms}.
 \end{equation}
 To show that $R_{F,G}$ does not vanish of infinite order at $x'=0$
 it is thus sufficient to show that $R_{\hat{F},\hat{G}}$ does not
 vanish identically.
 Therefore the theorem is proved if we can choose the $\alpha_m$ and
 $\beta_m$ in such a way that $R_{\hat{F},\hat{G}}$ does not
 vanish identically.
 
 Assume the resultant $R_{\hat{F},\hat{G}}$ vanishes for all
 linear combinations $\hat{F}$ and $\hat{G}$.
 Then for any $x'\in\R^{n-1}$ there exists a common root $\xi(x')$ of 
 the polynomials $w_1,\ldots,w_r$.
 By (\ref{ffw}) and (\ref{rf}) we have
 \begin{equation}
 \sum_{j=0}^{k}\frac{1}{j!}D_1^jy_0(x')\xi(x')^j = 0.
 \label{crf} 
 \end{equation}
 
 There is a nonempty open subset of $\R^{n-1}$ on which
 $\xi$ can be chosen such that it depends smoothly on $x'$.
 On this subset we apply $D_1$ to (\ref{crf}) and use (\ref{DPR}) to obtain
 \begin{equation}
 \sum_{j=0}^{k-1}\frac{\xi(x')^j}{j!}\cdot(1+\nabla\xi)\cdot
 D_1^{j+1}y_0(x') = 0.
 \label{dcrf} 
 \end{equation}
 
 If $v\in V$ is a vector in a finite dimensional Euclidean
 vector space $V$, then the element $1+v$ is invertible in the
 Clifford algebra $Cl(V)$ with inverse $\frac{1-v}{1+|v|^2}$.
 Thus (\ref{dcrf}) gives 
 \begin{equation}
 \sum_{j=0}^{k-1}\frac{\xi(x')^j}{j!}D_1^{j+1}y_0(x') = 0.
 \label{scrf} 
 \end{equation}
 Repeating this argument inductively we eventually get
 $$
 D_1^ky_0 = 0.
 $$
 But this means $y_k(x') = (v_1(0),\ldots,v_r(0)) = 0$, a contradiction.
 $\Box$

 \ep
 
 Mathematisches Institut
 
 Universit\"at Freiburg
 
 Eckerstr.\ 1
 
 79104 Freiburg
 
 Germany
 
 \vspace{0.3cm}
 
 {\em e-mail:}
 {\typewriter baer@mathematik.uni-freiburg.de}
 
 \end{document}